\documentclass[letterpaper, 10 pt, conference]{mys}  

\IEEEoverridecommandlockouts                              
\usepackage{graphics} 
\usepackage{epsfig} 
\usepackage{amsmath} 
\usepackage{amssymb}  
\usepackage[latin1]{inputenc}
\usepackage{color}

\newcommand{\eq}{\begin{equation}}
\newcommand{\eeq}{\end{equation}}
\newcommand{\eqn}{\begin{eqnarray}}
\newcommand{\eeqn}{\end{eqnarray}}

\newcommand{\bsea}{\begin{subeqnarray}}
\newcommand{\esea}{\end{subeqnarray}}
\newcommand{\nn}{\nonumber}

\newcommand{\Argmin}[1]{\,\underset{#1}{\mathrm{argmin}}\,}

\newcommand{\Ac}{ \mathcal{A}}
\newcommand{\Bc}{ \mathcal{B}}

\newcommand{\Lc}{ \mathcal{L}}

\newcommand{\Sc}{ \mathcal{S}}

\def\qed{\hfill \vrule height 7pt width 7pt depth 0pt \smallskip}

\newcounter{pippo}

\newtheorem{remark}{Remark}[section]
\newtheorem{teor}{Theorem}[section]
\newtheorem{corr}{Corollary}[section]
\newtheorem{propo}{Proposition}[section]
\newtheorem{lemm}{Lemma}[section]
\newtheorem{exam}{Example}
\newtheorem{obss}{Observation}
\newtheorem{probl}[pippo]{Problem}
\newtheorem{defn}{Definition}[section]
\newcommand{\teo}{\begin{teor}}
\newcommand{\eteo}{\end{teor}}
\newcommand{\cor}{\begin{corr}}
\newcommand{\ecor}{\end{corr}}
\newcommand{\prop}{\begin{propo}}
\newcommand{\eprop}{\end{propo}}
\newcommand{\lem}{\begin{lemm}}
\newcommand{\elem}{\end{lemm}}
\newcommand{\ex}{\begin{exam}}
\newcommand{\eex}{\end{exam}}
\newcommand{\pb}{\begin{probl}}
\newcommand{\epb}{\end{probl}}
\newcommand{\df}{\begin{defn}}
\newcommand{\edf}{\end{defn}}
\newcommand{\aprop}{\begin{apropo}}
\newcommand{\eaprop}{\end{apropo}}
\newcommand{\alem}{\begin{alemm}}
\newcommand{\ealem}{\end{alemm}}
\newcommand{\rem}{\begin{remark}}
\newcommand{\erem}{\end{remark}}
\newcommand{\oss}{\begin{obss}}
\newcommand{\eoss}{\end{obss}}

\title{\LARGE \bf A New Recursive Least-Squares Method\\ 
 with Multiple Forgetting Schemes}

\author{Francesco Fraccaroli, Andrea Peruffo and Mattia Zorzi
\thanks{This work has been partially supported by the FIRB project ``Learning
meets time'' (RBFR12M3AC) funded by MIUR.}
\thanks{F. Fraccaroli is with the Dipartimento di Ingegneria dell'Informazione, Universit\`a degli studi di
Padova, via Gradenigo 6/B, 35131 Padova, Italy
        {\tt\small francesco.fraccaroli.2@studenti.unipd.it}}%
\thanks{A. Peruffo is with the Dipartimento di Ingegneria dell'Informazione, Universit\`a degli studi di
Padova, via Gradenigo 6/B, 35131 Padova, Italy
        {\tt\small andrea.peruffo@studenti.unipd.it}}%
\thanks{M. Zorzi is with the Dipartimento di Ingegneria dell'Informazione, Universit\`a degli studi di
Padova, via Gradenigo 6/B, 35131 Padova, Italy
        {\tt\small zorzimat@dei.unipd.it}}%
}

\begin{document}

\maketitle
\thispagestyle{empty}
\pagestyle{empty}

\begin{abstract}
We propose a recursive least-squares method with multiple forgetting schemes to track time-varying model parameters which change with different rates. Our approach hinges on the reformulation of the classic recursive least-squares with forgetting scheme as a regularized least squares problem. A simulation study shows the effectiveness of the proposed method.
\end{abstract}

\section{Introduction}

Recursive identification methods are essential in system identification, 
\cite{LJUNG_SODERSTROM_1983,YOUNG_2011,SOLO_KONG_ADAPTATIVE_1995,ljung1981analysis,SODERSTROM_ML_1972,LJUNG_2002}. In particular, they are able to track variations of the model parameters over the time. This task is fundamental in adaptive control, \cite{aastrom2013adaptive,ljung1990adaptation,widrow1985adaptive}.

Recursive least-squares (RLS) methods with forgetting scheme represent a natural way to cope with recursive identification. These approaches can be understood as a weighted least-squares problem wherein the old measurements are exponentially discounted through a parameter called forgetting factor. Moreover, in \cite{CAMPI_1994} their tracking capability has been analysed in a rigorous way.

In this paper, we deal with models having time-varying parameters which change  with different rates. Many applications can be placed in this framework. An example is the automation of heavy duty
vehicles, \cite{VAHIDI_2005}. In this problem, it is required to estimate the vehicle mass and the road grade.
The former is almost constant over the time, whereas the latter is time-varying. Other examples are the control of strip temperature for heating furnace, \cite{yoshitani1998model}, and the self-tuning cruise control, \cite{oda1991practical}.

In those applications the RLS with forgetting scheme provides poor performances. A refinement of this method is the RLS with directional forgetting scheme, \cite{HAGGLUND_1985,Kulhavy_1987,Bittanti_1990,cao_schwartz_1999}. 
Roughly speaking, such approach fixes the problem that the incoming information is not uniformly distributed over all parameters. However, this nonuniformity
is not equivalent to the presence of parameters with different changing rates, \cite{VAHIDI_2005}. Indeed, it is 
possible to construct models with parameters having different changing rates and with incoming information uniformly distributed over all parameters. Thus, also RLS with directional forgetting scheme provides poor  performances.

An {\em ad-hoc} remedy to estimate parameters with different changing rates is the RLS with vector-type forgetting (or selective forgetting) scheme, \cite{saelid_1983,saelid_1985,PARKUN_1990,PARKUM_1992}. The idea of the above method is to introduce many forgetting factors reflecting the different rates of the change of the parameters. Finally, an {\em ad-hoc} modification of the above method has been presented in \cite{VAHIDI_2005}.

In this paper, we propose a new RLS with multiple forgetting schemes. Our method is based on the reformulation 
of the classic RLS with forgetting scheme as a regularized least-squares problem. It turns out that the current parameters vector minimizes the current prediction error 
plus a penalty term. The latter is the weighted distance between the current and the previous value of the parameters vector. Moreover, the weight matrix is updated at each time step and the updating law depends on the forgetting factor. This simple observation leads us to generalize this updating to multiple forgetting factors reflecting the different changing  rates of the parameters. 
Moreover, we provide three updating laws drawing inspiration on machine learning. For simplicity we will consider SISO models because the extension to MIMO ones is straightforward. Finally, simulation show the effectiveness of our method.

The remainder of the content in the paper is organized as follows. In Section~\ref{sec:review}, we present the state of the art about RLS with forgetting scheme and with vector-type forgetting scheme. The reformulation of the RLS and the three different updating laws are explained in Section \ref{sec:RLS}. The performance comparisons between these methods are illustrated in Section \ref{sec:simulation}. Conclusions are drawn in Section \ref{sec:conclusion}.

\section{State of the art}\label{sec:review}

Consider a SISO linear, discrete time, time-varying, system
\begin{equation} \label{eq:model}
\Ac_{t}(z^{-1}) y(t) =  \Bc_{t}(z^{-1}) u(t) + e(t),
\end{equation}
where $e(t)$ is additive noise with variance $\sigma^2$ and $u(t)$ is a stationary Gaussian process independent of $e(t)$.

$\Ac_t(z^{-1})$ and $\Bc_t(z^{-1})$ are time-varying polynomials whose degrees are $n$ and $m$ respectively:
\begin{subequations}
\begin{align}
& \Ac_t(z^{-1}) = 1 + \sum_{i=1}^n a_{t, i} z^{-i},  \\
& \Bc_t(z^{-1}) = \sum_{i=1}^m b_{t, i} z^{-i},
\end{align}
\end{subequations}
where $z$ is the shift operator.

Assume to collect the data \eq Z^N:=\{ y(1), u(1) \dots y(N), u(N) \}.\eeq We would estimate $\Ac_t(z^{-1})$ and $\Bc_t(z^{-1})$ at each time step $t$ given $Z^{t}$.  We define
\begin{equation}
\theta_t = [ -a_{t, 1} \ldots - a_{t, n} \;b_{t, 1} \ldots b_{t, m} ]^T
\end{equation}
as the vector containing the parameters of $\Ac_t(z^{-1})$ and $\Bc_t(z^{-1})$. Let $\Phi_t$ denote the regression matrix
\begin{equation}
\Phi_t =
\begin{bmatrix}
\varphi(t)^T \\
\vdots \\
\varphi(max(n, m) + 1)^T
\end{bmatrix}
\end{equation}
where $\varphi(t) = [y(t-1) \ldots y(t-n)\; u(t-1) \ldots u(t-m)]^T$.

Let $\mathbf{y}_t$ be the vector of observations
\begin{equation}
\mathbf{y}_t = [ y(t)  \dots y(max(n,m) + 1)]^T \quad 
\end{equation}
and in similar way $\mathbf{e}_t$ be the noise vector
\begin{equation}
\mathbf{e}_t = [e(t)  \dots e(max(n,m) + 1)]^T. \quad 
\end{equation}
A common way to solve such a problem relies on the RLS with forgetting scheme, \cite{LJUNG_SODERSTROM_1983,YOUNG_2011}, where $\hat \theta_t$
is given by
\begin{equation}
\label{eq:class-rpem}
\hat{\theta}_t = \Argmin{\theta_t} V(\theta_t,t),
\end{equation}
and the loss-function is \eq \label{loss_function}V(\theta,t)=\sum_{s=1}^t \lambda^{t-s} (y(s)-\varphi(s)^T \theta).\eeq
Here, the forgetting factor $\lambda\in[0,1]$ operates as an exponential weight which decreases for the more remote data.

Problem (\ref{eq:class-rpem}) admits the recursive solution
\begin{subequations}
\begin{align}
& R_t = \lambda R_{t-1} + \varphi(t) \varphi(t)^T, \label{eq:Rtrec} \\
& \hat{\theta}_t = \hat{\theta}_{t-1} + R_t^{-1} \varphi(t)( y(t) - \varphi(t)^T \hat{\theta}_{t-1} ).
\end{align}
\end{subequations}
Moreover, if we define $P_t = R_t^{-1}$ we obtain the equivalent recursion
\begin{subequations}
\begin{align}
& \hat{\theta}_t = \hat{\theta}_{t-1} + K_t ( y(t) - \varphi(t)^T \hat{\theta}_{t-1} ), \\
& K_t = \frac{P_{t-1}\varphi(t)}{\lambda + \varphi(t)^T P_{t-1} \varphi(t)}, \\
& P_t = \frac{1}{\lambda} (I - K_t \varphi(t)^T)P_{t-1}.
\end{align}
\end{subequations}

In the case that the parameters in ARX model (\ref{eq:model}) vary with a different rate it is desirable to assign different forgetting factors. The RLS with vector-type forgetting scheme, \cite{saelid_1985,PARKUN_1990}, consists of scaling $P_t$ by a diagonal matrix $\Lambda$ of forgetting factors
\begin{equation}
P_t = \Lambda^{- \frac{1}{2}} (I - K_t \varphi(t)^T)P_{t-1} \Lambda^{-\frac{1}{2}}
\end{equation}
where $\Lambda = diag(\lambda_1 \dots \lambda_p)$ with $p = n + m$. Therefore, $\lambda_i$ is the forgetting factor reflecting the changing  rate of the $i$-th parameter. Finally, an {\em ad-hoc} modification of the update law for the gain $K_t$ of the RLS has been proposed in \cite{VAHIDI_2005}. In this case the parameters to estimate are two. Such method conceptually separates the error due to the parameters in two parts in the objective function (\ref{loss_function}), that is one part contains the error due to the parameter with faster changing  rate and the second one the error due to the parameter with slower changing  rate. Then two different forgetting factors have been applied for each term.

\section{RLS with multiple forgetting schemes} \label{sec:RLS}
In this Section, we introduce our RLS for model whose parameters have different changing  rates.
Our approach hinges on the following observation.

\prop \label{lemma} Problem (\ref{eq:class-rpem}) is equivalent to the following problem:
\begin{align} \label{eq:reform}
\hat{\theta}_t = & \Argmin{\theta_t} (y(t)-\varphi(t)^T \theta_t )^2 + \nonumber \\
& + \lambda (\theta_t-\hat{\theta}_{t-1})^T R_{t-1} (\theta_t-\hat{\theta}_{t-1}),
\end{align}
with updating law (\ref{eq:Rtrec}). \eprop The proof is given in Appendix \ref{sec:AppA}.

Proposition \ref{lemma} shows that the RLS with forgetting scheme can be understood as regularized least squares problem. More precisely, the first term in the objective function minimizes the prediction error at time $t$, whereas the penalty term minimizes the distance between $\theta_t$ and the previous estimate $\hat \theta_{t-1}$ according to the weight matrix $\lambda R_{t-1}$. Moreover, the weight matrix is updated according to the law (\ref{eq:Rtrec}).

It is then natural to allow a more general structure for the weight matrix $\lambda R_{t-1}$ and its updating law (\ref{eq:Rtrec}). Let $F_{\lambda}(\cdot)$ be the forgetting map defined as follows
\begin{align*}
F_{\lambda} : & \ \Sc_p^{+} \ \rightarrow \ \Sc_p^{+} \\
& R_{t-1} \mapsto F_{\lambda}(R_{t-1}),
\end{align*}
where $\Sc_p^+$ denotes the cone of positive definite matrices of dimension $p$ and $\lambda = [\lambda_1 \dots \lambda_p]^T \in \mathbb{R}^p$ is the forgetting vector
with $0 < \lambda_i < 1$ $i=1\ldots p$ forgetting factor of the $i$-th parameter.

Therefore, given $\hat \theta_{t-1}$, we propose the following estimation scheme for $\theta_t$
\begin{subequations}
\begin{align}
& \hat{\theta_t} = \Argmin{\theta_t} (y(t)-\varphi(t)^T \theta_t )^2 + \nonumber \\
& \qquad + (\theta_t-\hat{\theta}_{t-1})^T F_{\lambda}(R_{t-1}) (\theta_t-\hat{\theta}_{t-1}),\label{eq:reformF} \\
& R_t = F_{\lambda}(R_{t-1}) + \varphi(t) \varphi(t)^T.  \label{eq:RtF}
\end{align}
\end{subequations}

\prop \label{propo} The solution to (\ref{eq:reformF}) with updating law (\ref{eq:RtF}) admits the recursive solution
\begin{subequations}
\begin{align}
& \hat{\theta}_t = \hat{\theta}_{t-1} + K_t (y(t) - \varphi(t)^T\hat{\theta}_{t-1})\label{eq:thetaF}, \\
& K_t = R_t^{-1} \varphi(t)  \label{eq:KtF}, \\
& R_t =F_\lambda(R_{t-1})+\varphi(t)\varphi(t)^T  \label{eq:PtF}.
\end{align}
\end{subequations}
Moroever,  $K_t$ can be updated in the equivalent way:
\begin{subequations}
\begin{align}
& K_t = \frac{F_{\lambda}(P_{t-1}^{-1})^{-1}\varphi(t)}{1 + \varphi(t)^T F_{\lambda}(P_{t-1}^{-1})^{-1}\varphi(t)}  \label{eq:KtF2}, \\
& P_t = (I - K_t \varphi(t)^T )F_{\lambda}(P_{t-1}^{-1})^{-1} \label{eq:PtF2}
\end{align}
\end{subequations} where $P_t=R_t^{-1}$.
\eprop

The proof is given in Appendix \ref{sec:AppB}.

To design the forgetting map $F_\lambda$ we consider the following result whose proof can be found in \cite{RASMUSSEN_WILLIAMNS_2006}.
\prop Consider $A,B\in\mathcal{S}_p^+$. Let $C$ be a symmetric matrix of dimension $p$ such that
\eq  [C]_{ij}=[A]_{ij}[B]_{ij},\;\; i,j=1\ldots p.\eeq Then, $C\in\mathcal{S}_p^+$.
\eprop

In view of the above result, a natural structure for $F_\lambda$ would be
\eq [F_\lambda(R_{t-1})]_{ij}=[R_{t-1}]_{ij} [Q_\lambda]_{ij}\eeq
where $Q_\lambda\in\Sc_p^+$. Note that, $Q_\lambda$ can be understood as a kernel
matrix with hyperparameters $\lambda$ in the context of machine learning, \cite{RASMUSSEN_WILLIAMNS_2006,WAHBA1990}.
Next, we design three types of maps drawing inspiration on the diagonal kernel, the tuned/correlated kernel, \cite{EST_TF_REVISITED_2012},
and the cubic spline kernel, \cite{ WAHBA1990}.

\subsection{Diagonal updating} \label{sec:diag_updating}
Consider the ARX model (\ref{eq:model}) with $m=1$ and $n=1$, therefore we only have two parameters.
Let $\theta_{t,1}$ and $\theta_{t,2}$ denote the parameter of $\Ac_t(z^{-1})$  
and $\Bc_t(z^{-1})$, respectively. Moreover, the vector containing the two parameters is defined as 
$\theta_t=\left[
\begin{array}{cc}
\theta_{t,1}& \theta_{t,2}\\
\end{array}
\right]^T
$. We assume that the changing  rate of  $\theta_{t,1}$ is slow over the interval $[1,N]$, whereas the changing  rate of $\theta_{t,2}$ is faster.
The simplest idea is to decouple the parameters in the penalty term in (\ref{eq:reformF}). We associate the forgetting factor $\lambda_1$ to $\theta_{t,1}$ and $\lambda_2$ to $\theta_{t,2}$ with $\lambda_1>\lambda_2$. Let
\begin{equation}
R_{t-1} =
\begin{bmatrix}
R_{t-1, 1} & R_{t-1, 12} \\
R_{t-1, 12} & R_{t-1, 2}
\end{bmatrix}.
\end{equation}
Then, if we define
\begin{equation}
F_{\lambda,DI}(R_{t-1}) =
\begin{bmatrix}
\lambda_1 R_{t-1, 1} & 0 \\
0 & \lambda_2 R_{t-1, 2}
\end{bmatrix}
\end{equation}
the penalty term in (\ref{eq:reformF}) becomes
\eq \lambda_1(\theta_{t,1}-\hat \theta_{t-1,1})^2 R_{t-1,1}
+\lambda_2(\theta_{t,2}-\hat \theta_{t-1,2})^2R_{t-1,2}
\eeq that is the parameters of $\Ac_t(z^{-1})$ and the ones of $\Bc_t(z^{-1})$ have been decoupled in the penalty term.

This simple example leads us to consider the diagonal updating
\begin{equation*}
[F_{\lambda, DI} (R_{t-1})]_{i,j}= \left\{ \begin{matrix}
0 & if \ \lambda_i \neq \lambda_j \\
[R_{t-1}]_{i,j} \lambda_i &  \hbox{otherwise}
\end{matrix} \right. .
\end{equation*}
Finally, it is worth noting that in the special case that $p=2$ we obtain the method proposed in \cite[\textit{formulae} (22) and (23)]{VAHIDI_2005}.

\subsection{Tuned/Correlated updating} \label{sec:TC} 
We consider again the example of Section  \ref{sec:diag_updating}.
The changing  rate of $R_{t-1,12}$ depends on the changing  rates of $\theta_{t,1}$
and $\theta_{t,2}$. Hence, it is reasonable to forget past values of $R_{t-1,12}$
with the fastest changing rate between the one of $\theta_{t,1}$ and $\theta_{t,2}$. Therefore,
we weigh $R_{t-1,12}$ with the forgetting factor $\lambda_2$
\begin{equation*}
F_{\lambda, TC} (R_{t-1}) =
\begin{bmatrix}
\lambda_1 R_{t-1, 1} & \lambda_2 R_{t-1, 12} \\
 \lambda_2 R_{t-1, 12} & \lambda_2 R_{t-1, 2}
\end{bmatrix}.
\end{equation*}
Moreover, the corresponding penalty term is
\eqn & \lambda_1(\theta_{t,1}-\hat \theta_{t-1,1})^2 R_{t-1,1}
 +\lambda_2(\theta_{t,2}-\hat \theta_{t-1,2})^2R_{t-1,2}\nn\\
& + 2\lambda_2(\theta_{t,1}-\hat \theta_{t-1,1}) (\theta_{t,2}-\hat \theta_{t-1,2}) R_{t-1,12}.
\eeqn Thus, the weight of the cross term is dominated by the smallest forgetting factor.
Therefore, in the general case, a reasonable updating law is: 
\begin{equation}
[F_{\lambda, TC} (R_{t-1})]_{i, j} = \min(\lambda_i, \lambda_j) [R_{t-1}]_{i, j}.
\end{equation}

\subsection{Cubic Spline updating} \label{sec:CS}
Consider the example of Section \ref{sec:diag_updating}.
We want to construct an updating such that the weight of the cross term in the penalty term (\ref{eq:reformF}) is not totally dominated by the forgetting factor $\lambda_2$. More precisely, we want that the weight of the cross term is also influenced by $\lambda_1$. We consider $Q_\lambda$ as a cubic spline like kernel matrix
\eq \label{eq:Q_CS}[Q_\lambda]_{ij}=\min \left[ \ \frac{l_i^2}{2} \left(l_j - \frac{l_i}{3}\right), \frac{l_j^2}{2} \left( l_i - \frac{l_j}{3}\right) \ \right] \eeq
where $l_1,l_2> 0$, $i=1,2$, is a function of $i$ to be determined. In our case we want that
\eq [Q_\lambda]_{ii}=\frac{l_i^3}{3}\eeq
is equal to $\lambda_i$ for $i=1,2$. Therefore, we obtain
\eq \label{eq:li}l_i=\sqrt[3]{3\lambda_i},\;\; i=1,2.\eeq
In this way, we built a forgetting map whose cross term is penalized by a blend of $\lambda_1$ and $\lambda_2$.

\begin{figure}[h]
	\centering
		\includegraphics[width=\columnwidth]{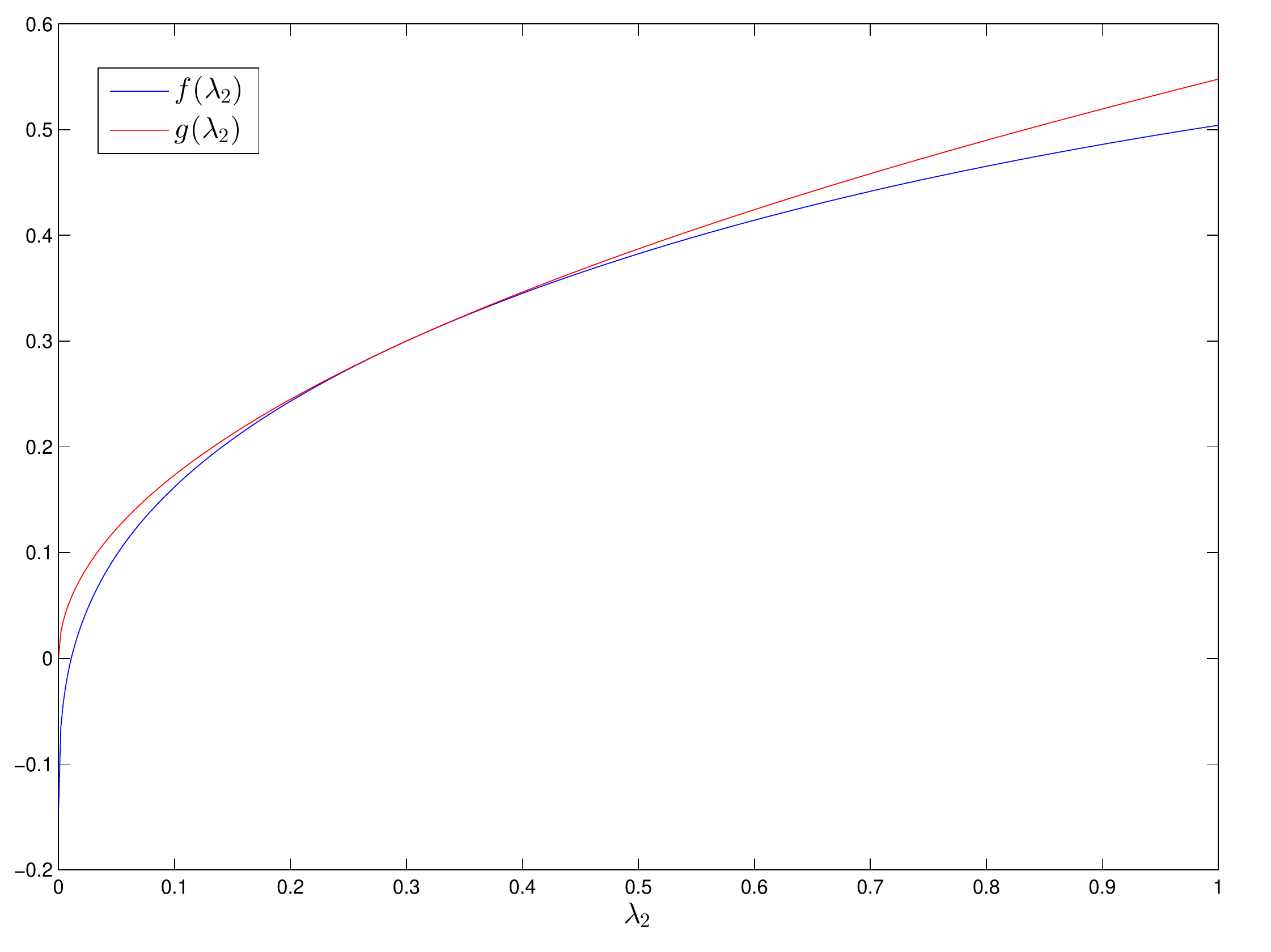}
		\caption{\textit{Comparison between $f(\cdot)$ and $g(\cdot)$}}
		\label{fig:cs}		
\end{figure}

\rem One could also consider the matrix $\tilde Q_\lambda$ such that $[\tilde{Q}_\lambda]_{ji}=\sqrt{\lambda_i\lambda_j}$, $i,j=1\ldots p$. To compare (\ref{eq:Q_CS}) and $\tilde Q_\lambda$ assume that $\lambda_1$ is fixed equal to 0.3, whereas $\lambda_2$ can vary over the interval $[0,1]$. In Figure \ref{fig:cs} we depict the functions $f(\lambda_2) = \frac{l_1^2}{2} \left( l_2- \frac{l_1}{3} \right)$ and $g(\lambda_2) = \sqrt{\lambda_1 \lambda_2}$. As one can see $f(\cdot)$  takes smaller values than the ones of $g(\cdot)$ , that is the influence of the smallest forgetting factor is more marked in $f(\cdot)$.

Thus, by plots evidence, (\ref{eq:Q_CS}) provides a blend of $\lambda_1$ and $\lambda_2$ in which the
influence of $\lambda_2$ (forgetting factor associated to the parameter with the fastest changing  rate) is more marked than the one in $\tilde Q_\lambda$.\erem

In the general case, therefore the updating law becomes 
\begin{align}
[F_{\lambda, CS} (R_{t-1})]_{i, j} := & [R_{t-1}]_{i, j} \times \nonumber \\
& \times \min \left[ \ \frac{l_i^2}{2} \left(l_j - \frac{l_i}{3}\right), \frac{l_j^2}{2} \left( l_i - \frac{l_j}{3}\right) \ \right].
\end{align}
where $l_i = \,\sqrt[3]{3\lambda_i}$, $i=1\ldots p$.

\section{Simulations Results} \label{sec:simulation}

In this section we analyse the performance of the RLS with multiple forgetting schemes that we presented in Section \ref{sec:RLS}. The experiment has been performed using MATLAB as the numerical platform.

\subsection{Data generation}\label{sec:gen_data}

We consider a discrete-time, time-varying ARX model described in (\ref{eq:model}), with $n = 2$ and $m = 2$.
Here, the parameters in $\Ac_t(z^{-1})$ vary faster than the ones in $\Bc_t(z^{-1})$. To this aim, nine stable polynomials $\Ac^{(j)}(z^{-1})$, $j=1\ldots 9$, and two stable polynomials $\Bc^{(k)}(z^{-1})$, $k=1,2$, have been defined. We considered the time interval $[1,N]$ with $N=160$. The polynomial $\Bc_t(z^{-1})$ is generated as a smooth time varying convex combination of $\Bc^{(1)}(z^{-1})$ and $\Bc^{(2)}(z^{-1})$. Regarding $\Ac_t(z^{-1})$, we split the interval $[1,N]$ in eight sub-interval and at the $j$-th interval $\Ac_t(z^{-1})$ is generated as a smooth time varying convex combination of $\Ac^{(j)}(z^{-1})$ and $\Ac^{(j+1)}(z^{-1})$.

Finally, the input $u(t)$ is generated as a realization of white Gaussian noise with unit variance and filtered with a 10$^{th}$ order Butterworth low-pass filter. Starting from random initial condition, the output $y(t)$ is collected and corrupted by an additive white Gaussian noise with variance $\sigma^2=0.01$.

\subsection{Proposed Methods}

The method we consider are: 

\begin{itemize}
  \item RARX: this is the classic RARX algorithm implemented in \verb#rarx.m# in the MATLAB System identification Toolbox, \cite{LJUNG_SYS_ID_1999};
  \item VF: this is the RLS with vector-type forgetting scheme described at the end of Section \ref{sec:review};
  \item DI: this is the RLS algorithm with diagonal updating of Section \ref{sec:diag_updating};
  \item TC: this is the RLS algorithm with tuned/correlated updating of Section \ref{sec:TC};
  \item CS: this is the RLS algorithm with cubic spline updating of Section \ref{sec:CS}.
\end{itemize}
For each method $m=2$ and $n=2$, that is the estimated ARX models have the same order of the true one. Regarding VF, DI, TC and CS
we set
\eq \lambda=\left[
              \begin{array}{cccc}
                \lambda_1 & \lambda_1 & \lambda_2 & \lambda_2 \\
              \end{array}
            \right]^T
\eeq that is $\lambda_1$ is the forgetting factor for the parameters in $\Ac_t(z^{-1})$ and
 $\lambda_2$ is the forgetting factor for the parameters in $\Bc_t(z^{-1})$.

\subsection{Experiment setup}

We consider a study of $500$ runs.
For each run, we generate the data as described in Section \ref{sec:gen_data} and we compute $\hat \theta_t$ with the five methods. More precisely, for each method (VF, DI, TC and CS) we compute $\hat \theta_t$ for twenty values of $\lambda_1$ and $\lambda_2$ uniformly sampled over the interval $[0.1,1]$. Then, we pick $\lambda_1^\circ$ and $\lambda_2^\circ$ which maximize the {\em one step ahead coefficient of determination} (in percentage)
\begin{equation}
 \mathrm{COD}= \left( 1 - \frac{ \frac{1}{N} \sum_{t = 1}^N (y(t) - \hat{y}(t) )^2}{ \frac{1}{N} \sum_{t = 1}^N (y(t) - \overline{y}_N )^2 } \right) \times 100
\end{equation}
where $\hat{y}(t)$ is the predicted value of $y(t)$ based on the ARX model with $\Ac_{t-1}(z^{-1})$ and $\Bc_{t-1}(z^{-1})$, and $\overline{y}_N$ is the sample mean of the output data. It is worth noting that the
performance index $\mathrm{COD}$ is used for time invariant models. On the other hand, it provides a rough idea whether the estimated model is good or not and it allows to choose reasonable values for $\lambda_1^\circ$ and $\lambda_2^\circ$.
Then, for $\lambda^\circ_1$ and $\lambda_2^\circ$ we compute the corresponding  {\em average track fit} (in percentage)
\begin{equation}
\mathrm{ ATF }  = \left( 1- \frac{1}{N}\sum_{t=1}^N \frac{|| {\hat{\theta}}_t - \mathbf{\theta}_t||}{||{\theta}_t||} \right) \times 100.
\end{equation}
Regarding RARX, we use the procedure above with one forgetting factor.

\subsection{Results}

\begin{figure*}[th!]
	\centering
		\includegraphics[width=2\columnwidth]{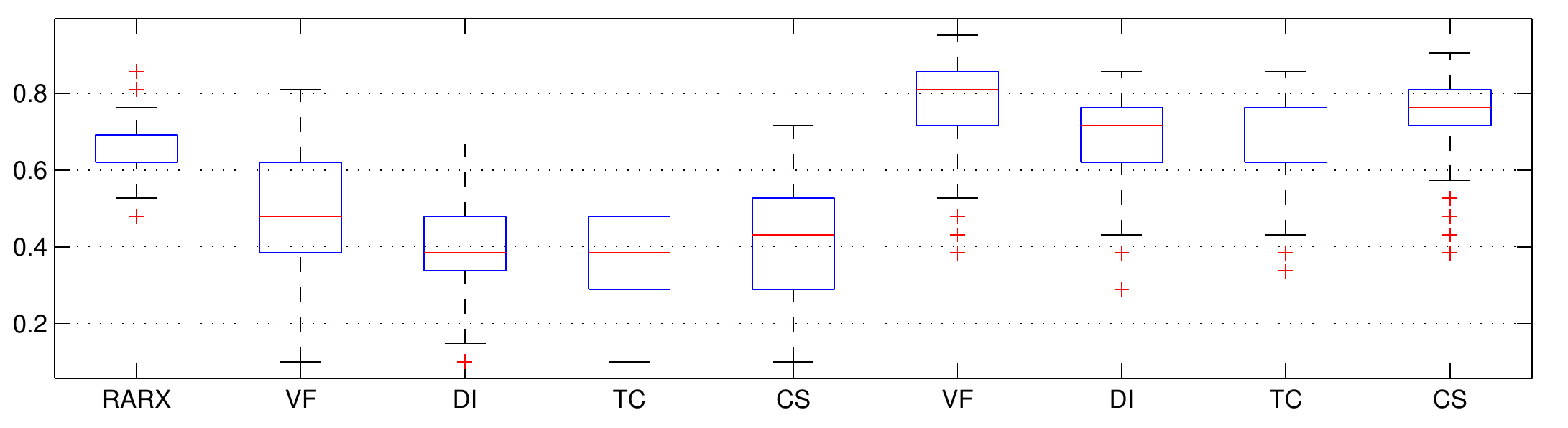}
		\caption{\textit{Forgetting factors of the different algorithms. First column: forgetting factor of RARX. Second-fifth column: forgetting factor $\lambda_1$ for VF, DI, TC and CS. Sixth-last column: forgetting factor $\lambda_2$ for VF, DI, TC and CS.}}
		\label{fig:lambda}		
\end{figure*}

In Figure \ref{fig:lambda} are shown the values of $\lambda$. The first boxplot refers to the values chosen by the classic RARX algorithm, from the second to the fifth the values of the forgetting factor $\lambda_1$ referring to the parameters of $\Ac_t(z^{-1})$ are represented, while the last ones refer to the forgetting factor $\lambda_2$ related to the parameters of $\Bc_t(z^{-1})$. Since the parameters of $\Ac_t(z^{-1})$ varies faster than the ones of $\Bc_t(z^{-1})$, its forgetting factors are smaller than the respective others, as expected. On the other hand, the classic RARX has not the possibility to choose different forgetting factors so its best choice is to take an intermediate value among the ones picked by the proposed algorithms.

\begin{figure}[h]
	\centering
		\includegraphics[width=\columnwidth]{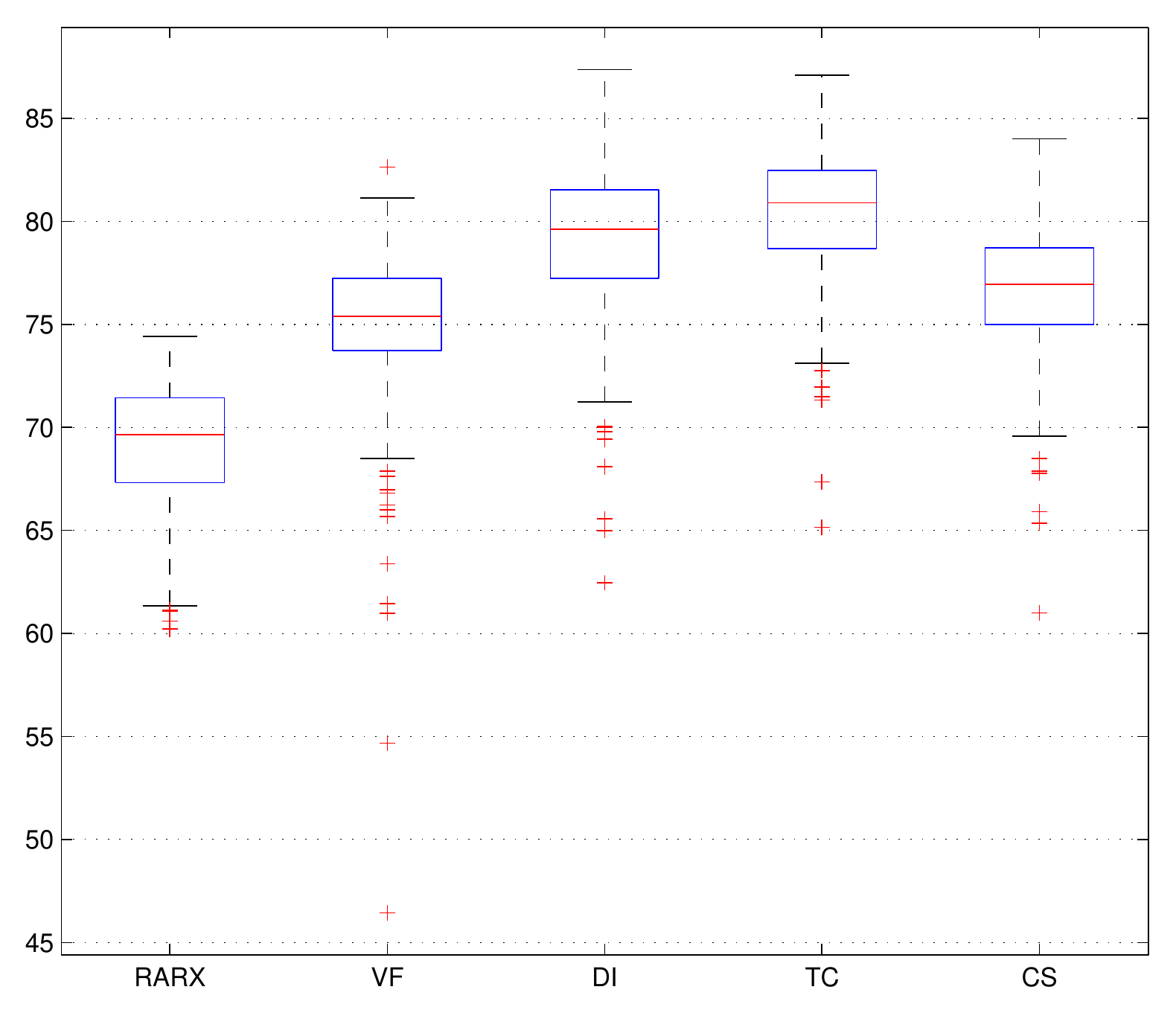}
		\caption{\textit{Average track fit of the different algorithms.}}
		\label{fig:apf}		
\end{figure}

In Figure \ref{fig:apf} are depicted the average track fit indexes. All the proposed algorithms have better performances than RARX and VF, anyway it is possible to highlight that the TC updating shows the best results. This fact suggests that the most efficient weight for the cross terms in the penalty term in (\ref{eq:reformF})  is the smallest forgetting factor between the eligible ones, as occurs in the TC algorithm.

\begin{figure}[h]
	\centering
		\includegraphics[width=\columnwidth]{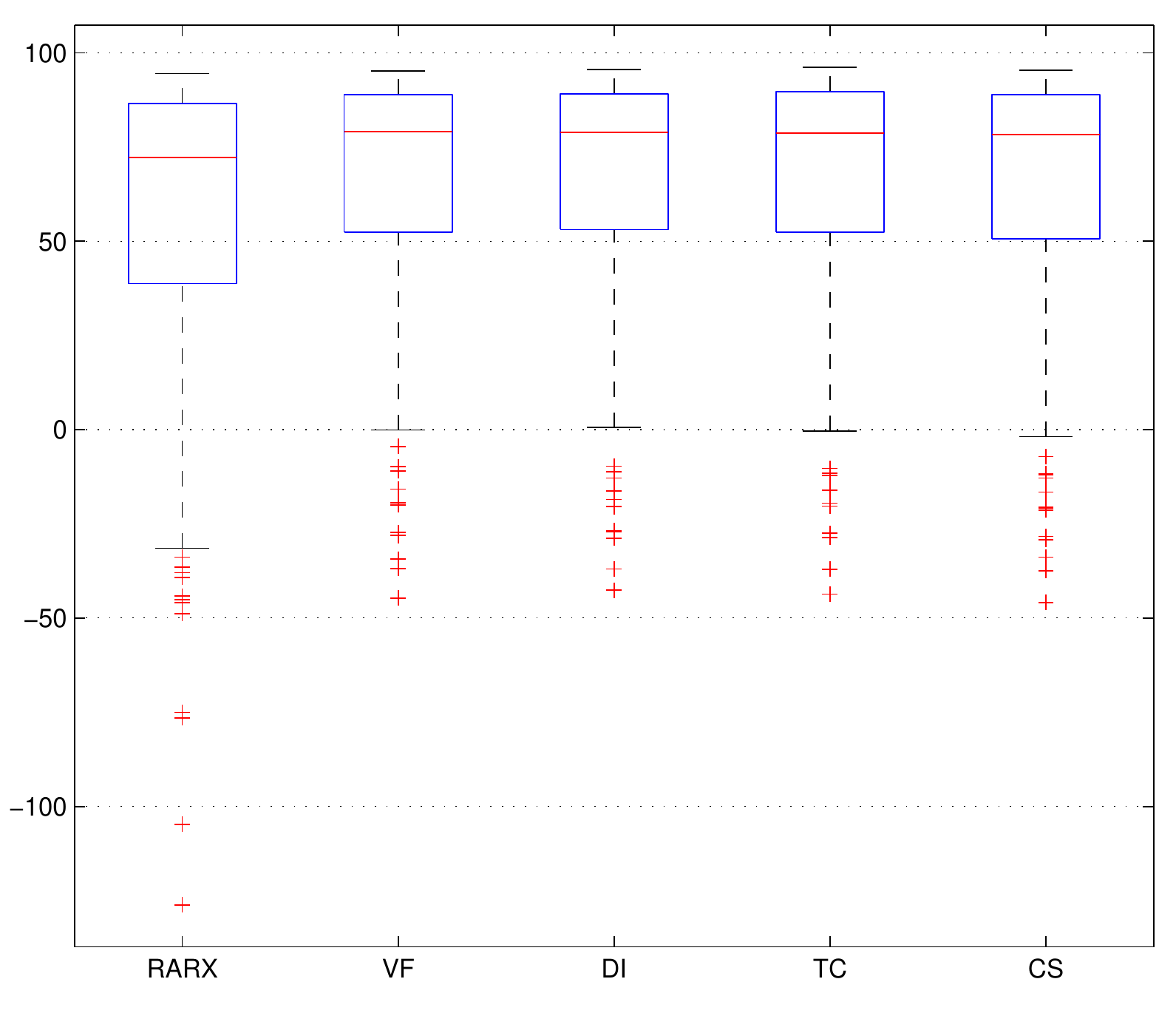}
		\caption{\textit{One step ahead coefficient of determination of the different algorithms.}}
		\label{fig:fit}		
\end{figure}

Figure \ref{fig:fit} illustrates the $\mathrm{COD}$ indexes. Once again the proposed algorithms outperforms the classic RARX method: if we focus on the average value of the boxplots the difference is around 5\%. In terms of outliers we can underline that RARX reaches $-100 \%$ in the worst case scenario, while the proposed methods never go below $-55 \%$.

\section{Conclusions} \label{sec:conclusion}

We presented a reformulation of the classic RLS algorithm, which can be split into the minimization of the current prediction error and the minimization of a quadratic function which penalizes the distance between the current and the previous value of the estimate. This reformulation is strictly connected to an updating equation which provides the weight matrix of the quadratic function: to change the updating equation given by the classic algorithm means to substitute the map that connects the present weight matrix to the past one. This permits to model multiple forgetting factors to improve the estimation of parameters with different changing  rates.

In this paper we provide three different updating laws. Simulations show that these algorithms outperforms the conventional ones thanks to the proposed updating law which allows the presence of several forgetting factors. Therefore, multiple forgetting factors  seem to be the key to a more efficient identification. It is worth noting that the challenging step is the choice of such forgetting factors. Therefore, the next research direction will concern the estimation of such parameters from the collected data.

\appendix

\section{Appendix A}

\subsection{Proof of Proposition \ref{lemma}} \label{sec:AppA}
Let $
Q_t = diag(1 \dots \lambda^{t-1}).
$ Consider 

\begin{align*}
\hat{\theta}_t & = \Argmin{\theta_t} \sum_{i=1}^t (y(i) - \varphi(i)^T \theta_t)^2 \lambda^{t-i} \\
& = \Argmin{\theta_t} (y(t) - \varphi(t)^T \theta_t)^2 + \lambda \sum_{i=1}^{t-1} (y(i) - \varphi(i)^T \theta_t )^2 \lambda^{t-i-1} \\
& = \Argmin{\theta_t}(y(t) - \varphi(t)^T \theta_t)^2 + \\
& \quad + \lambda \sum_{i=1}^{t-1} [ y(i) - \varphi(i)^T \theta_t  + \varphi(i)^T\hat{\theta}_{t-1} - \varphi(i)^T \hat{\theta}_{t-1} ]^2 \lambda^{t-i-1} \\
& = \Argmin{\theta_t} (y(t) - \varphi(t)^T \theta_t)^2 + \lambda \sum_{i=1}^{t-1} [ \ (\varphi(i)^T(\theta_t - \hat{\theta}_{t-1})) ^2  + \\
& \quad - 2( y(i) - \varphi(i)^T \hat{\theta}_{t-1})\varphi^T(i) (\theta_t- \hat{\theta}_{t-1}  ) \ ] \lambda^{t-i-1},
\end{align*}
where the term $(y(i) - \varphi(i)^T \hat{\theta}_{t-1})^2$ has been omitted because it does not depend on $\theta_t$.

The last equation can be rewritten as 
\begin{align}
\hat{\theta}_t & = \Argmin{\theta_t}  (y(t)- \varphi(t)^T \theta_t)^2 + \nonumber \\ 
& + \lambda [ \ || \theta_t - \hat{\theta}_{t-1} ||^2_{\Phi_{t-1}^T Q_{t-1} \Phi_{t-1} } + \nonumber  \\
& \quad -2 (\theta_t - \hat{\theta}_{t-1})^T \Phi_{t-1}^T Q_{t-1} (\mathbf{y}_{t-1} - \Phi_{t-1} \hat{\theta}_{t-1}) \ ].
\label{eq:thm2prod}
\end{align}

It is not difficult to see that 
\begin{equation*}
\hat \theta_{t-1} = \Argmin{\theta_{t-1}} || \mathbf{y}_{t-1} - \Phi_{t-1} \theta_{t-1} ||^2_{ Q_{t-1}  }.
\end{equation*}

Therefore, it must hold the following equation (by optimality condition)
\begin{equation*}
\Phi_{t-1}^T Q_{t-1} (\mathbf{y}_{t-1} - \Phi_{t-1} \hat{\theta}_{t-1}) = 0,
\end{equation*}
so (\ref{eq:thm2prod}) becomes 
\begin{align*}
\hat{\theta}_t  & = \Argmin{\theta_t}  (y(t)- \varphi(t)^T \theta_t)^2 + \nonumber \\
& \quad + \lambda || \theta_t - \hat{\theta}_{t-1} ||^2_{\Phi_{t-1}^T Q_{t-1} \Phi_{t-1} } .
\end{align*}

Finally, it is sufficient to observe that $R_t = \Phi_t^T Q_t \Phi_t$. \qed\\

\subsection{Proof of Proposition \ref{propo}} \label{sec:AppB}

The objective function in (\ref{eq:reformF}) is
\begin{equation*}
\Lc(\theta_t) := (y(t)- \varphi(t)^T \theta_t)^2 +  || \theta_t - \hat{\theta}_{t-1} ||^2_{F_{\lambda}(R_{t-1}) }.
\end{equation*}
Therefore, the optimal solution takes the form
\begin{align*}
\hat{\theta}_t & = [ \varphi(t) \varphi(t)^T + F_{\lambda}(R_{t-1}) ]^{-1} (\varphi(t) y(t) + F_{\lambda}(R_{t-1}) \hat{\theta}_{t-1}). \end{align*}
Finally, (\ref{eq:thetaF}), (\ref{eq:KtF}), (\ref{eq:PtF}), (\ref{eq:KtF2}) and (\ref{eq:PtF2})
can be derived along similar lines used for the classic RLS with forgetting scheme.\qed\\

\bibliographystyle{plain}

\begin{thebibliography}{10}

\bibitem{aastrom2013adaptive}
K.~{\AA}str{\"o}m and B.~Wittenmark.
\newblock {\em Adaptive control}.
\newblock Courier Corporation, 2013.

\bibitem{Bittanti_1990}
S.~Bittanti, P.~Bolzern, and M.~Campi.
\newblock Convergence and exponential convergence of identification algorithms
  with directional forgetting factor.
\newblock {\em Automatica}, 26(5):929--932, 1990.

\bibitem{CAMPI_1994}
M.~Campi.
\newblock Performance of rls identification algorithms with forgetting factor:
  A phi-mixing approach.
\newblock {\em Journal of Mathematical Systems, Estimation and Control},
  4(3):1--25, 1994.

\bibitem{cao_schwartz_1999}
L.~Cao and H.~Schwartz.
\newblock A novel recursive algorithm for directional forgetting.
\newblock In {\em Proceedings of the American Control Conference}, volume~2,
  pages 1334--1338, 1999.

\bibitem{EST_TF_REVISITED_2012}
T.~Chen, H.~Ohlsson, and L.~Ljung.
\newblock On the estimation of transfer functions, regularizations and gaussian
  processes-revisited.
\newblock {\em Automatica}, 48(8):1525--1535, 2012.

\bibitem{HAGGLUND_1985}
T.~H\"{a}gglund.
\newblock Recursive estimation of slowly time-varying parameters.
\newblock In {\em Proc. IFAC Symposium on Identification and System Parameter
  Estimation}, pages 1137--1142, York, 1985.

\bibitem{SOLO_KONG_ADAPTATIVE_1995}
V.~Kong and X.~Solo.
\newblock {\em Adaptive signal processing algrithms}.
\newblock Prentice Hall, New Jersey, 1995.

\bibitem{Kulhavy_1987}
R.~Kulhav{\`y}.
\newblock Restricted exponential forgetting in real-time identification.
\newblock {\em Automatica}, 23(5):589--600, 1987.

\bibitem{ljung1981analysis}
L.~Ljung.
\newblock Analysis of a general recursive prediction error identification
  algorithm.
\newblock {\em Automatica}, 27(1):89--100, 1981.

\bibitem{LJUNG_SYS_ID_1999}
L.~Ljung, editor.
\newblock {\em System Identification (2Nd Ed.): Theory for the User}.
\newblock Prentice Hall, New Jersey, 1999.

\bibitem{LJUNG_2002}
L.~Ljung.
\newblock Recursive identification algorithms.
\newblock {\em Circuits, Systems and Signal Processing}, 21(1):57--68, 2002.

\bibitem{ljung1990adaptation}
L.~Ljung and S.~Gunnarsson.
\newblock Adaptation and tracking in system identification--{A} survey.
\newblock {\em Automatica}, 26(1):7--21, 1990.

\bibitem{LJUNG_SODERSTROM_1983}
L.~Ljung and T.~S{\"o}derstr{\"o}m.
\newblock {\em Theory and Practice of Recursive Identification}.
\newblock MIT Press, 1983.

\bibitem{oda1991practical}
K.~Oda, H.~Takeuchi, M.~Tsujii, and M.~Ohba.
\newblock Practical estimator for self-tuning automotive cruise control.
\newblock In {\em American Control Conference}, pages 2066--2071, 1991.

\bibitem{PARKUN_1990}
J.~Parkum, N.~Poulsen, and J.~Holst.
\newblock Selective forgetting in adaptive procedures.
\newblock In {\em The 11th IFAC World Congress in Tallinn}, pages 180--185,
  1990.

\bibitem{PARKUM_1992}
J.~Parkum, N.~Poulsen, and J.~Holst.
\newblock Recursive forgetting algorithms.
\newblock {\em International Journal of Control}, 55(1):109--128, 1992.

\bibitem{RASMUSSEN_WILLIAMNS_2006}
C.~Rasmussen and C.~Williams.
\newblock {\em {Gaussian Processes for Machine Learning}}.
\newblock MIT Press, 2006.

\bibitem{saelid_1985}
S.~Saelid, O.~Egeland, and B.~Foss.
\newblock A solution to the blow-up problem in adaptive controllers.
\newblock {\em Modeling, Identification and Control}, 6(1):1--36, 1985.

\bibitem{saelid_1983}
S.~Saelid and B.~Foss.
\newblock Adaptive controllers with a vector variable forgetting factor.
\newblock In {\em The 22nd IEEE Conference on Decision and Control}, pages
  1488--1494, 1983.

\bibitem{SODERSTROM_ML_1972}
T.~Söderström.
\newblock An on-line algorithm for approximate maximum likelihood
  identification of linear dynamic systems.
\newblock Technical Report 7308, Department of automatic control, Lund
  Institute of technology, Lund, Sweden, 1973.

\bibitem{VAHIDI_2005}
A.~Vahidi, A.~Stefanopoulou, and H.~Peng.
\newblock Recursive least squares with forgetting for online estimation of
  vehicle mass and road grade: theory and experiments.
\newblock {\em Vehicle System Dynamics}, 43(1):31--55, 2005.

\bibitem{WAHBA1990}
G.~Wahba.
\newblock {\em Spline Models for Observational Data.}
\newblock Society for Industrial and Applied Mathematics, 1990.

\bibitem{widrow1985adaptive}
B.~Widrow and S.~Stearns.
\newblock {\em Adaptive signal processing}.
\newblock Prentice-Hall, New Jersey, 1985.

\bibitem{yoshitani1998model}
N.~Yoshitani and A.~Hasegawa.
\newblock Model-based control of strip temperature for the heating furnace in
  continuous annealing.
\newblock {\em IEEE Transactions on Control Systems Technology}, 6(2):146--156,
  1998.

\bibitem{YOUNG_2011}
P.~Young.
\newblock {\em Recursive estimation and time-series analysis: An introduction
  for the student and practitioner}.
\newblock Springer Science \& Business Media, 2011.

\end{thebibliography}

\end{document}